\pdfoutput=1
\documentclass{ifacconf}
\usepackage[T1]{fontenc}
\usepackage[protrusion=true,expansion=true,final]{microtype} 
\usepackage{graphicx,color,enumerate,makecell}
\usepackage{natbib}
\usepackage{textcomp}

\usepackage{amsmath,amssymb,amsfonts,ae,bm,mathtools}  

\usepackage[algo2e,ruled]{algorithm2e}

\renewcommand*{\thefootnote}{\fnsymbol{footnote}}
\DeclareMathAlphabet{\mathdutchcal}{U}{dutchcal}{m}{n}

\usepackage{float}
\usepackage{ctable}

\newtheorem{proposition}{Proposition}

\theoremstyle{definition}

\theoremstyle{plain}

\newcommand*{\tran}{{\mkern-1.5mu\mathsf{T}}}

\begin{document}
% --------------------------------------------------------------------
% *************************** FRONT MATTER ***************************
% --------------------------------------------------------------------
\begin{frontmatter}
% -----------------

\title{Block particle filters for state estimation of stochastic reaction-diffusion systems} 

 \thanks[footnoteinfo]{``\textcopyright 2023 the authors. This work has been accepted to IFAC for publication under a Creative Commons Licence CC-BY-NC-ND''}
% goes here. Paper titles should be written in uppercase and lowercase
% letters, not all uppercase.}  

\author[First]{Jos\'e Augusto F. Magalhães}  
\author[First]{Otac\'ilio B. L. Neto}  
\author[First]{Francesco Corona}

\address[First]{School of Chemical Engineering, Aalto University, Finland. \\(e-mail: \{augusto.magalhaes, otacilio.neto, francesco.corona\}@aalto.fi)}

\begin{abstract}	
	In this work, we consider a differential description of the evolution of the state of a reaction-diffusion system under environmental fluctuations. We are interested in estimating the state of the system when only partial observations are available. To describe how observations and states are related, we combine multiplicative noise-driven dynamics with an observation model. More specifically, we ensure that the observations are subjected to error in the form of additive noise. We focus on the state estimation of a Belousov-Zhabotinskii chemical reaction. We simulate a reaction conducted in a quasi-two-dimensional physical domain, such as on the surface of a Petri dish. We aim at reconstructing the emerging chemical patterns based on noisy spectral observations. For this task, we consider a finite difference representation on the spatial domain, where nodes are chosen according to observation sites. We approximate the solution to this state estimation problem with the Block particle filter, a sequential Monte Carlo method capable of addressing the associated high-dimensionality of this state-space representation.
\end{abstract}

\end{frontmatter}
%===============================================================================

%===============================================================================
\section{Introduction}
%===============================================================================
% Reaction-diffusion systems refers to processes of chemical reactions between spatially distributed substances whose molecules flow according to their concentration gradient. This class of systems have been of great interest to the scientific community as they form an essential basis for studying self-assembly and self-organization phenomena occurring in nature. 

% React-diffusion
Since the introduction of the ``Turing pattern'' concept, reaction-diffusion models have played an important role in theoretical biology \citep{Turing1952}. Such representations have been used to explain several naturally-occurring pattern formations: developing structures in embryos \citep{Kondo2010,Green2015}, skin and pigmentation patterns \citep{Kondo2009,Nakamasu2009}, and tumour cell growth \citep{Ferreira2002,Hogea2008}. Additionally, reaction-diffusion models have also been employed to describe wildfire spread dynamics \citep{Asensio2002} and to synthesize textures in computer graphics \citep{Witkin1991}. 
Moreover, if the state of a reaction-diffusion system can be determined in real-time, one could use feedback theory to control the pattern-formation process for several purposes \citep{Ardizzone2013}. 
Adopting these models in real biochemical applications has been challenging: such models must rely on partial observations of the true system through indirect methods (e.g. light spectrography), and their typical realisations are high-dimensional stochastic systems. Consequently, the task of estimating the state is bound to incur a heavy computational burden, if at all computable. 

% Classical methods 
%In order to overcome this limitation, methods such as the Ensemble Kalman Filter (EnKF, \citep{Evensen2003}) have been developed and successfully applied in many fields. Specifically, this method has been popular for applications in which partial differential equation (PDE) models are discretised into a large number of ordinary differential equations (ODE) to be simulated in a computer \citep{Anderson2009,Aanonsen2009}. Additionally, the EnKF algorithm can be easily extended to a parallel formulation which makes it efficient for implementation in real systems \citep{Elias2019}.

% Particle filters   
Particle filters  (PFs) are popularly used for state estimation of general non-linear stochastic dynamical systems \citep{Chopin2020}. An ensemble of particles equipped with importance weights undergoes a sampling procedure every time a new observation is available \citep{Kitagawa1996}. Due to the recursive nature of this particle representation, the ensemble is bound to present a degeneracy problem \citep{Musso2001}.
Some techniques to mitigate this phenomenon include introducing diversity among particles \citep{Gilks2001} or resampling a section of their paths over a fixed lag \citep{Doucet2006}. However, even when equipped with these techniques, classical particle filters applied to models of high-dimensional state space will suffer from degeneracy even after a few iterations \citep{VanLeeuwen2003}. An alternative approach consists on leveraging the decay-of-correlations property of many high-dimensional systems by performing full-state estimation based on lower-dimensional partitions of the state-space.
% Another strategy to reduce degeneracy is to work with local low-dimensional portions of the model. 
This class of methods, denoted local particle filters, include the Multiple particle filter (MPF, \cite{Djuric2007}) and, more recently, the Block particle filter (BPF, \cite{Rebeschini2015}).

% Work summary
Herein, we aim at estimating the concentrations of chemical substances in general reactive-diffusion systems using noisy spectral observations. We focus on state-space representations arising from finite-difference approximations of the spatio-temporal dynamics. We address this high-dimensional problem by using BPFs based on disjoint partitions of the spatial domain. We illustrate the performance of the estimator  for the task of reconstructing the emerging chemical patterns of a Belousov-Zhabotinskii-type reaction-diffusion system, the Oregonator \citep{Field1974}. From our results, the \emph{vanilla} BPF is unable to reconstruct the state for relatively small ensemble sizes. However, we show that the BPF is able to provide accurate estimates of the state under a specific choice of sampling procedure.
%it is clear that the choice of sampling procedure is crucial to reconstruct the state accurately when using a relatively small ensemble of particles. Under a specific choice of sampling procedure, the BPF yields accurate results for a problem where classical  particle filters would have been degenerate. 

% Flow
The paper is organised as follows: Section~\ref{sec: Reaction_Diffusion} overviews reaction-diffusion systems and their representation using finite difference methods. Section~\ref{sec: Preliminaries} discusses PFs for general and high-dimensional state spaces. Section~\ref{sec: Results} presents the simulation results of applying the state estimator on a benchmark reaction-diffusion model. Finally, Section~\ref{sec: Concluding remarks} lists our final remarks and some future research directions. % when applying the BPF under this model and

%===============================================================================

% ===============================================================================
\section{Reaction-diffusion systems} \label{sec: Reaction_Diffusion}
% ===============================================================================

We consider sets of chemical species $\mathcal{S} = \{\mathcal{S}_1,\dots,\mathcal{S}_{N_{\mathcal{S}}}\}$ distributed over some bounded space $U \subseteq \mathbb{R}^3$ which interact according to chemical reaction networks of the form
\begin{equation} \label{eq: CRN}
    \sum_{n_s = 1}^{N_{\mathcal{S}}} \underline{S}_{n_r, n_s} \mathcal{S}_{n_s} \xrightarrow{\kappa_{n_r}} \sum_{n_s = 1}^{N_{\mathcal{S}}} \overline{S}_{n_r, n_s} \mathcal{S}_{n_s},  \quad n_r = 1,\dots, N_R,
\end{equation}
with rate constants $\kappa \in \mathbb{R}^{N_R}_{\geq 0}$, and matrices $\underline{S} \in \mathbb{N}^{N_R\times N_{\mathcal{S}}}$ and $\overline{S} \in \mathbb{N}^{N_R\times N_{\mathcal{S}}}$ of stoichiometric coefficients for the reactants and products, respectively. 
From a system-analytical perspective, this process has dynamics described in terms of concentrations $z(u,t) = ([\mathcal{S}_1](u,t),\dots,[\mathcal{S}_{N_{\mathcal{S}}}](u,t))\in \mathbb{R}_{\geq 0}^{N_{\mathcal{S}}}$ at every space coordinates $u \in U$ and time-instant $t \in \mathbb{R}_{\geq 0}$. The evolution of concentrations can then be described as
\begin{equation} \label{eq: Reaction_Diffusion}
    \partial^{t} z(u,t) = \left[ (\overline{S}-\underline{S})^{\tran} \nu\big(z(u,t)\big) + D_z \nabla^2 z(u,t) \right] \partial t,
\end{equation}
given reaction rates $\nu(\cdot) = (\nu_1(\cdot), \ldots, \nu_{N_R}(\cdot)) \in \mathbb{R}^{N_R}$, matrix of diffusion coefficients $D_z \in \mathbb{R}^{N_{\mathcal{S}}\times N_{\mathcal{S}}}$, and Laplacian $\nabla^2 z(u,t) = \sum_{i=1}^{3} \partial z^2 / \partial u_i^2$. For each $n_r$-th reaction in network Eq. \eqref{eq: CRN}, the function $\nu_{n_r}$ ($n_r = 1,\ldots,N_R$) follows the Law of Mass Action \citep{Murray2002},
\begin{equation}
    \nu_{n_r} \hspace*{-0.1em}\big( z(u,t) \big) = \kappa_{n_r} \displaystyle\prod_{n_s=1}^{N_{\mathcal{S}}} z_{n_s} \hspace*{-0.1em}(u,t)^{\underline{S}_{n_r, n_s}}.
\end{equation} 
Finally, we consider homogeneous Neumann conditions $(\partial z / \partial u)\vert_{u = \overline{u}} = 0$ for all boundary points $\overline{u} \in \partial U$.

The concentrations of each substance are assumed to be only observed indirectly through their mixture. Specifically, we assume a spectroscopy process emitting
\begin{equation} \label{eq: Y_Process}
    y(u,t|\lambda) = \displaystyle\sum_{n_s=1}^{N_{\mathcal{S}}} \phi_{\mathcal{S}_{n_s}}\hspace*{-0.2em}(\lambda) z_{n_s}\hspace*{-0.1em} (u,t), 
\end{equation}
given wavelength $\lambda \in \mathbb{R}_{\geq 0}$ and functions $\phi_{\mathcal{S}_{n_s}} : \mathbb{R}_{\geq 0} \to \mathbb{R}$ describing the characteristic response spectrum associated with each $n_s$-th species ($n_s = 1,\ldots,N_{\mathcal{S}}$). This corresponds to light absorbance being linearly related to the concentration of substances as in Beer's Law \citep{Boaz2005}. Hereafter, we consider the response $\Phi_{\mathcal{S}_{n_s}}\hspace*{-0.2em} (\Lambda) = (\phi_{\mathcal{S}_{n_s}}(\lambda_1),\ldots,\phi_{\mathcal{S}_{n_s}}(\lambda_{N_{\Lambda}}))$ for specific wavelengths $\Lambda = (\lambda_1,\ldots,\lambda_{N_{\Lambda}})$, so that the measurement process is expressed by $y(u,t) = H_z z(u,t)$ with $H_z = [\Phi_{\mathcal{S}_1} \cdots \Phi_{\mathcal{S}_{N_{\mathcal{S}}}}]$. 

The reaction-diffusion system Eq. \eqref{eq: CRN} with measurement process Eq. \eqref{eq: Y_Process} can be represented by the state-space model 
\begin{subequations}  
\begin{align} 
    \partial^t z(u,t) &= f_{\theta_x}^z(z(u,t), \nabla^2 z(u,t)) \partial t; \label{eq: CRN_SS_PDE_A}\\
               y(u,t) &= H_{\theta_y}^z z(u,t),            \label{eq: CRN_SS_PDE_B}
\end{align}\label{eq: CRN_SS_PDE}%
\end{subequations}
with state $z : U \times \mathbb{R}_{\geq 0} \to \mathbb{R}_{\geq 0}^{N_{\mathcal{S}}}$ and measurements $y : U \times \mathbb{R}_{\geq 0}\to \mathbb{R}_{\geq 0}^{N_{\Lambda}}$. Function $f^z(\cdot)$ and matrix $H^z \in \mathbb{R}^{N_{\Lambda} \times N_{\mathcal{S}}}$ are respectively determined by fixed parameters $\theta_x$ and $\theta_y$ collected from Eqs. \eqref{eq: CRN}--\eqref{eq: Y_Process}.
 
In general, obtaining an exact solution for state-equation Eq. \eqref{eq: CRN_SS_PDE_A} is unpractical. Towards a numerical approach for its integration, we discretise the space $U$ into a lattice graph representing evenly spaced coordinates. State-dynamics can then be approximated by a collection of ordinary differential equations by the finite-difference method \citep{Strauss2007}. In the following, we overview this approach. 

\subsection{Finite-difference approximation}  \label{sec: Reaction_Diffusion_FDM}

We will restrict ourselves to quasi-two-dimensional spaces $U = \{ u \in \mathbb{R}^2 : 0 \leq u_{i=1,2} < \overline{U}\}$ with finite sizes $\overline{U} > 0$: these correspond to chemical mixtures in squared surfaces (e.g. a petri dish). Given some length $\Delta u > 0$, we discretise $U$ onto the square lattice $V = \{ v \in \mathbb{N}^2 : 1 \leq v_{i=1,2} \leq \overline{V}\}$ of evenly-spaced grid points with size $\overline{V} = (\overline{U}/\Delta u)$.
The spatially-discretised dynamics are then approximated as
\begin{equation} \label{eq: Reaction_Diffusion_Discretised}
    dz^{(v)} = \left[ (\overline{S}-\underline{S})^{\tran} \nu\big( z^{(v)} \big) + D_z \widetilde{\nabla}^2 z^{(v)} \right] dt, \quad \forall v \in V,
\end{equation}
with $z^{(v)}(t) = z\big((v_1,v_2)\Delta u,t\big)$, and the first-order approximation of the Laplacian $\nabla^2 z$ at $v = (v_1,v_2)$,
\begin{equation} \label{eq: 5Stencil_FDM}
    \widetilde{\nabla}^2 z^{(v)} = 
        \displaystyle\sum_{v'\in \mathcal{N}(v)} \frac{1}{\Delta u^2}\left( z^{(v')} - z^{(v)} \right),
        % \frac{z^{(v_1,v_2+1)}+z^{(v_1,v_2-1)}+z^{(v_1-1,v_2)}+z^{(v_1+1,v_2)}-4z^{(v_1,v_2)}}{\Delta u^2},
\end{equation}
where $\mathcal{N}(v) = \{ v' \in \mathbb{N}^2 : |v_1-v'_1| + |v_2-v'_2| = 1 \}$, the $1$-step neighbourhood of $v$.
The output equation becomes $y^{(v)}(t) = H_z z^{(v)}(t)$, with $y^{(v)}(t) = y((v_1,v_2)\Delta u, t)$. This measurement process thus corresponds to a typical imaging spectrography with \emph{pixels} assigned to grid points on $V$.

Finally, the reaction-diffusion model Eq. \eqref{eq: CRN_SS_PDE} is approximated by the conventional state-space representation
\begin{subequations}
\begin{align}
    dx(t) &= f_{\theta_x}(x(t)) dt; \label{eq: CRN_SS_ODE_A}\\
     y(t) &= H_{\theta_y} x(t),    \label{eq: CRN_SS_ODE_B}
\end{align}\label{eq: CRN_SS_ODE}%
\end{subequations}
with state-vector $x(t) = (z^{(v_1,v_2)}(t))_{v_1,v_2=1}^{\overline{V}} \in \mathbb{R}^{N_x}$ and output-vector $y(t) = (y^{(v_1,v_2)}(t))_{v_1,v_2=1}^{\overline{V}} \in \mathbb{R}^{N_y}$, obtained by collecting $(z^{(v)},y^{(v)})$ at every point in the square lattice $V$. The output matrix is $H_{\theta_y} = [I_{\overline{V}^2}\otimes\Phi_{\mathcal{S}_1} \cdots I_{\overline{V}^2}\otimes\Phi_{\mathcal{S}_{N_{\mathcal{S}}}}]$, where $M_1 \otimes M_2$ denotes the Kronecker product between matrices $M_1$ and $M_2$. The model thus corresponds to $N_x = N_{\mathcal{S}}\overline{V}^2$ state- and $N_y = N_{\Lambda}\overline{V}^2$ output-variables: A high-dimensional state space when $\overline{V}$ is large (i.e. a finely discretised space with $\Delta u \ll \overline{U}$).

% ===============================================================================

% ===============================================================================
\section{Preliminaries: Non-linear filtering} \label{sec: Preliminaries}
% ===============================================================================

We consider the stochastic differential equation (SDE) representation of an autonomous system,
\begin{subequations}
\begin{align}
    dX(t) &= f_{\theta_x}(X(t))dt + g_{\theta_b}(X(t)) dB(t),  \label{eq: SS_a}\\%
     Y(t) &= h_{\theta_y}(X(t)) + e(t), \label{eq: SS_b}%
\end{align} \label{eq: SS}%
\end{subequations}
with state equation Eq. \eqref{eq: SS_a} describing the time evolution of state process $X:[0,\infty)\times\Omega\to\mathbb{X}$ given its current value, stochastic driving process $B=(B(t))_{t\geq 0}$ and initial state $X(0) = \xi(\omega^*), \omega^*\in\Omega$.  
The output equation Eq. \eqref{eq: SS_b} describes how the current state is emitted as the output-vector $y(t) \in \mathbb{Y}$, after being corrupted by a measurement noise $e(t) \in \mathbb{R}^{N_y}$. 
We consider Polish spaces $\mathbb{X} = \mathbb{R}^{N_x}$ and $\mathbb{Y} = \mathbb{R}^{N_y}$. 
The parameter vector $\theta=(\theta_x, \theta_b, \theta_y)$ determines the functions $f_{\theta_x}(\cdot)$, $g_{\theta_b}(\cdot)$, and $h_{\theta_y}(\cdot)$. For the sake of exposition,  we omit the dependency on $\theta$ from here on.
For all $0\leq t'<t$, $B(t)-B(t')\sim \mathcal{N}(0,(t - t')I_{N_x})$ and $e(t) \sim \mathcal{N}(0,\Sigma_y)$ with known covariance $\Sigma_y$. 
Finally, we limit ourselves to linear equations $Y(t) = HX(t) + e(t)$ given matrix $H \in \mathbb{R}^{N_y \times N_x}$.

We are interested in estimating the states $(X(t))_{t \geq 0}$ based on a sequence of discrete-time measurements $(Y_1,\ldots,Y_k)$, with $Y_k = Y(k\Delta t)$ given a sampling time $\Delta t > 0$. Assuming a zero-order hold between measurements of the driving noise (that is, $g(X(t)) = g(X(t_{k-1}))$ for all $t \in [t_{k-1}, t_{k})$), the discrete-time state dynamics are represented by
\begin{equation} \label{eq: discretised-general-SDE}
    X_k = \underbrace{X_{k-1} + \int_{t_{k-1}}^{t_k}\hspace*{-0.5em} f(X(t)) dt }_{F(X_{k-1})} + g(X_{k-1})\Delta B_k, 
\end{equation}
where $X_k = X(k\Delta t)$, $\Delta B_k = \mathcal{N}(0,\Delta t \cdot I_{N_x})$, and $t_k = k\Delta t$. The discrete-time output equation is $Y_k = h(X_k) + e_k$ with measurement noise $e_k \sim \mathcal{N}(0,\Sigma_y)$.
The estimation problem is then formalised by the stochastic process $\pi=(\pi_k)_{k \in \mathbb{N}}$, where $\pi_{k} = \mathbb{P}(X_k \in \cdot | Y_1,\ldots, Y_k)$ is the filtering distribution summarising the uncertainty on $X_k$ given history $(Y_1,\ldots, Y_k)$.
The distributions $\pi_k$ are assumed to have density $p_k$ with respect to Lebesgue measure. 

We consider the Particle filter (PF, \cite{Chopin2020}) approach of computing Monte Carlo approximations $\pi^{N_p} = (\pi_k^{N_p})_{k \in \mathbb{N}}$ by recursively sampling $N_p$ candidates of $X_k$ based on observation $Y_k$ and previous $\pi^{N_p}_{k-1}$.
Accounting for high-dimensional spaces, we specialise this approach to a class of localised estimators known as the BPF \citep{Rebeschini2015}.
In the following, we overview these approaches.

% Subsections ====================================================================
% ===============================================================================
\subsection{Classical particle filters} \label{sec: Particle_Filter}
% ===============================================================================

We assume that $X$ is a Markov chain: for all $k\in\mathbb{N}$ and arbitrary Borel set $A\in\mathcal{B}(\mathbb{R}^{N_x})$, 
\begin{align}
\mathbb{P}(X_k\in A\vert X_{0:k-1}) = \mathbb{P}(X_k\in A\vert X_{k-1}).
\end{align}
We consider transition kernels $\text{P}_k : \mathbb{R}^{N_x}\times\mathcal{B}(\mathbb{R}^{N_x}) \to \mathbb{R}_{\geq 0}$ such that, for all $k\in \mathbb{N}$, $x\in\mathbb{R}^{N_x}$, and any $A\in\mathcal{B}(\mathbb{R}^{N_x})$,
\begin{align}
\text{P}_k(x, \text{A}) = \mathbb{P}(X_k\in\text{A}\vert X_{k-1}=x).\label{eq: transition-kernel}
\end{align}
The filter $\pi_k$ can thus be computed recursively by
\begin{align*}
	\pi_{k-1}\xrightarrow{\text{prediction}}\tilde{\pi}_{k} = \text{P}_k\pi_{k-1} \xrightarrow{\text{correction}} \pi_{k} = \text{C}_k\tilde{\pi}_{k},
\end{align*}
with the functions $\text{P}_k$ and $\text{C}_k$ satisfying
\begin{align}
		 	(\text{P}_k\pi_{k-1})(A) &\coloneqq \int \text{P}_k(x, A) \pi_{k-1}(dx), \label{eq: prediction}\\
	(\text{C}_k\tilde{\pi}_{k})(A) &\coloneqq \frac{\int_A\text{C}_k(x)\tilde{\pi}_{k}(dx)}{\tilde{\pi}_{k}(\text{C}_k)},\label{eq: correction}
\end{align}
for non-negative functions $\text{C}_k$ with $\tilde{\pi}_{k}(\text{C}_k)>0$. Here, $\tilde{\pi}_k(\text{C}) \coloneqq \int\text{C}(x)\tilde{\pi}_k(dx)$. 
In the general case, the recursion $\pi_k = \text{C}_k \text{P}_k \pi_{k-1}$ has no trivial solution. As such, the estimation problem can be solved by a Sequential Monte Carlo approach: The filter $\pi_k$ is approximated by an empirical distribution $\text{C}_k \text{S}_k^{N_p} \tilde{\pi}_{k}$ given the sampling operator $\text{S}_k^{N_p} \pi \coloneqq N_p^{-1} \sum_{n=1}^{N_p} \delta_{x^{n}}$ with samples $x^{1},\ldots,x^{N_p} \overset{i.i.d}{\sim} \pi$.

%The following exposition is adapted from \cite{Crisan2008}. 
Let particles $\{X^{(n)}_k\}_{n=1}^{N_p}$ be $N_p$ mutually independent stochastic processes, all independent of $Y_k$, solving Eq.~\eqref{eq: SS_a}. The pairs $(X^{(n)}_k, Y_k)$ are identically distributed with the same distribution as $(X_k,Y_k)$, $n=1,\ldots,N_p$. 
Now, let $\pi^{N_p} = (\pi_k^{N_p})_{k\in\mathbb{N}}$ be the sequence of empirical distributions
\begin{align} \label{eq: pf-approximation}
		\pi_k^{N_p} \,\overset{\Delta}{=} \, \frac{1}{N_p}\sum_{n=1}^{N_p} w^{(n)}_k \delta_{X^{(n)}_k},
\end{align}
where to each $X^{(n)}_k$ we assign a (normalised) weight $w^{(n)}_{k}$. 
Then, we have $\pi^{N_p} \to \pi$ almost surely at a rate slightly lower than $1/\sqrt{N_p}$ \citep{Crisan2008}.
As the filters $\pi_k$ are unavailable, we sample i.i.d. $X^{(1)}_k,\ldots,X^{(N_p)}_k$ from \emph{importance distributions} $\rho_k$ with densities of the form
\begin{align}
	q(x_{0:k}\vert y_{1:k}) = q(x_0)\prod_{k'=1}^{k}q( x_{k'}\vert y_{0:k'},  x_{0:k'-1}),
 	\label{eq: SIR-bayes}
\end{align} 
with $x_k$ and $y_k$ denoting a realisation of $X$ and $Y$ at time $t_k$, respectively.
In this setup, the unnormalised weights $\widetilde{w}^{(n)}(x_{0:k}) = w^{(n)}(x_{0:k-1}) \widetilde{w}^{(n)}_k$ ($n=1,\ldots,N_p$) can be computed recursively according to
\begin{align}
	\widetilde{w}^{(n)}_k = \frac{p(y_k\vert y_{1:k-1}, x^{(n)}_{0:k})p(x_k\vert x^{(n)}_{k-1})}{q(x_k\vert y_{1:k},x^{(n)}_{0:k-1})} \coloneqq l_k(x^{(n)}_{k},y_k),
	\label{eq: SIR-weights}
\end{align}  
whenever new observations $Y_k$ become available.
They are normalised through $w_k = [\sum_{n'=1}^{N_p} \widetilde{w}_k^{(n')}]^{-1}\widetilde{w}_k^{(n)}$.

The choice of $q(x_{0:k}\vert y_{1:k})$ is arbitrary: It only needs to have a support including that of $p(x_{0:k}\vert y_{0:k})$.  
Conventionally, we sample particles from the transition distribution of the dynamics Eq.~\eqref{eq: SS_a}, i.e. $q(x_k\vert y_{1:k},x_{0:k-1}) = p(x_k\vert x_{k-1})$.
The importance weights are then based on the likelihood, 
\begin{align}
	\widetilde{w}^{(n)}_k \propto w^{(n)}_{k-1} p(y_k\vert x^{(n)}_{k}).
	\label{eq: standard weights}
\end{align}
This choice of proposal density is conventional but not very informative. In practice, only a few selected particles might be assigned relevant weights. In the extreme situation, termed degeneracy, $w^{(n)}_k \approx 1$ for a single $n$, and $w^{(n')}_k \approx 0$ for all $n'\neq n$. The likelihood is then poorly approximated with large variance $\texttt{Var}_{\rho}(\widetilde{w})$. 

\subsubsection*{Optimal sampling distribution} \label{optimal-IS}
The degeneracy problem can be addressed by choosing an importance distribution that minimises the variance of the (unnormalised) weights conditioned not only on $x_{0:k}$ and $y_{0:k-1}$, but also on $y_k$. 

\begin{proposition}{\citep{Doucet-SMC-importance-sampling:20}.} \label{th: Optimal_IS_Distribution}
	The (optimal importance) distribution that minimises $\texttt{Var}_{\rho}(\widetilde{w})$ has density
	\begin{align}
		q(x_k\vert y_{1:k},x_{0:k-1}) =\frac{p(y_k\vert y_{1:k-1}, x_{0:k})p(x_k\vert x_{k-1})}{p(y_k\vert y_{1:k-1},x_{0:k-1})},
		\label{eq: SIR-proposal}
	\end{align}  
	where the associated unnormalised weights satisfy
	\begin{align}
		\widetilde{w}_k &\propto 
			w_{k-1}\underbrace{\int p(y_k \vert x_{k-1:k})p(x_k\vert x_{k-1})\,dx_k}_{p(y_k\vert x_{k-1})} .
		\label{eq: SIR-weights-optimal}
	\end{align}  
\end{proposition}

For processes $X_k = F(X_{k-1}) + g(X_{k-1})\Delta B_k$ with measurement models $Y_k = HX_k + e_k$ (given $e_k \sim \mathcal{N}(0,\Sigma_y)$), the optimal importance distribution has the analytical form
\begin{equation} \label{eq: Optimal_Importance_Distribution}
	X_k\vert Y_{k},X_{k-1} \sim \mathcal{N}(M_k^{\text{opt}}, \Sigma_k^{\text{opt}})
\end{equation}
with mean $M_k^{\text{opt}}$ and covariance $\Sigma_k^{\text{opt}}$ computed by
\begin{align}
	\Sigma_B 	 		  &= (g(x_{k-1})g(x_{k-1})^{\tran})\Delta t,\\
	\Sigma_k^{\text{opt}} &= (\Sigma_B^{-1} + H^{\tran} \Sigma_y^{-1} H)^{-1},\\
	M_k^{\text{opt}} 	  &= \Sigma_k^{\text{opt}}(\Sigma_B^{-1}F(x_{k-1})) + H^{\tran} \Sigma_y^{-1} y_k.
\end{align}
Moreover, the likelihood distribution satisfies 
\begin{align} \label{eq: optimal weights}
	Y_k \vert X_{k-1} \sim \mathcal{N}(HF(x_{k-1}),\ \Sigma_B + H\Sigma_yH^{\tran}).
\end{align} 
For positive systems, a common modelling assumption is that $g(X_{k-1}) = \Sigma_B^{1/2} =  \texttt{diag}(X_{k-1})\Sigma_x^{1/2}$, given a $\Sigma_x^{1/2} \succ 0$: This assures a unique and positive solution $X$ of Eq.~\eqref{eq: SS_a}, \citep{Yang2020}. As such, variances in Eqs.~\eqref{eq: Optimal_Importance_Distribution} and \eqref{eq: optimal weights} can be computed efficiently for high-dimensional systems when $H$ has some special structure (such as in Section \ref{sec: Reaction_Diffusion}).
% ===============================================================================
\subsection{Block particle filtering} \label{section: the block particle filter}
% ===============================================================================

In general, regardless of the choice of importance distribution, classical particle filters face degeneracy issues when applied to high-dimensional systems \citep{Snyder2015}. A scalable solution is to design filters that are spatially-localised: Dynamics and observations at
a spatial location are assumed to depend only on state-variables associated with its neighbourhood \citep{Rebeschini2015}. 

Consider the pair $(X_k, Y_k)$ at each time $t_k$ to be a random field $(X_k, Y_k)_{v\in V}$ indexed by a (finite) undirected graph $\mathcal{G} = (V,E)$.
The spaces $\mathbb{X}$ and $\mathbb{Y}$ of $X_k$ and $Y_k$ and the filter $\pi_k$ can be expressed in the product form
\begin{align}
	\mathbb{X} = {\textstyle\prod}_{v\in V}\mathbb{X}^v,\quad\mathbb{Y} = {\textstyle\prod}_{v\in V}\mathbb{Y}^v,\quad\pi_k = {\textstyle\bigotimes}_{v\in V}\pi_k^v.
\end{align}
Here, $\pi_k^v$ is the conditional distribution on state-space $\mathbb{X}^v$ and $\bigotimes_{v\in V}$ is the product (of measures) over a vertex set $V$.
The transition kernel $\text{P}_k$ in Eq.~\eqref{eq: transition-kernel} and likelihood function $l$ in Eq.~\eqref{eq: SIR-weights} are then decomposed in their localised forms
\begin{align*}
	\text{P}_k(x, A) = {\textstyle\prod}_{v\in V}\text{P}_k^v(x, A^{(v)}),\, l_k(x, y) = {\textstyle\prod}_{v\in V}l_k^v(x^{(v)}, y^{(v)}).
\end{align*}
We focus on the scenario where the graph $\mathcal{G}$ represents a square lattice (as in Section~\ref{sec: Reaction_Diffusion_FDM}) with vertex set 
\begin{align}
	V = \{ v \in \mathbb{N}^2 : 1 \leq v_{i=1,2} \leq \overline{V} \} \quad (\overline{V}\in\mathbb{N}).
\end{align}
Considering the set of non-overlapping blocks
\begin{align}
	\mathcal{K} = \{(v_0 + \{1,\ldots,\overline{V}_b\}^2) \cap V : v_0 \in \overline{V}_b\mathbb{N}^2\},
\end{align}
with size $\overline{V}_b < \overline{V}$, the vertices of $\mathcal{G}$ can be partitioned as
\begin{align*}
	V = \bigcup_{V_b \in \mathcal{K}} V_b, \quad V_b \cap V_{b'} = \emptyset \text{ for } V_b \neq V_{b'},\quad  V_b, V_{b'} \in \mathcal{K}.
\end{align*}
Under such decomposition, we define $X^{(V_b,n)} = \otimes_{v\in V_b} X^{(v,n)}$ and $Y^{(V_b)} = \otimes_{v\in V_b} Y^{(v)}$ to be the ($n$-th particle) state- and output-variables associated with block $V_b \in \mathcal{K}$.

Let the blocking operator be $\text{B}\pi_k \coloneqq \bigotimes_{V_b \in \mathcal{K}} \text{B}^{V_b}\pi_k$, where $\text{B}^{V_b}\pi_k$ denotes the marginal distribution of $\pi_k$ on $\prod_{v \in {V_b}}\mathbb{X}^{v}$. 
In this case, we compute filters $\pi_k$ recursively by
\begin{align*}
	\pi_{k-1} \xrightarrow{\text{prediction}} \tilde{\pi}_{k} = \text{P}_k \pi_{k-1} \xrightarrow{\text{correction}} \pi_{k} = \text{C}_k\text{B}\tilde{\pi}_{k}.
\end{align*}
The BPF algorithm is then implemented as in Algorithm~\ref{alg: BPF}\footnote[5]{In its standard formulation, Algorithm~\ref{alg: BPF} assumes $P_k^{v}$ and $l_k^{v}$ to be obtained from the appropriate marginals with densities $p(\tilde{x}_k \vert x_{k-1}^{(n)})$ and $p(y_k \vert \tilde{x}_k^{(n)})$, respectively. When the optimal importance distribution is used, the marginals are taken from Eqs.~\eqref{eq: Optimal_Importance_Distribution} and \eqref{eq: optimal weights}.}.
In this approach, the convergence properties of the particle filter become dependent on the cardinality of individual blocks, $|V_b| = \overline{V}_b^2$ ($V_b \in \mathcal{K}$), rather than on the cardinality of the original lattice, $|V| = \overline{V}^2$.

\begin{algorithm2e}
	\caption{Block particle filter}	\label{alg: BPF}
	\SetKwFunction{Resample}{Resample}
	Initialise $\pi_0^{N_p}$ with a desirable distribution;

	\For{$k = 1,2,\ldots,K$}{
		Resample $X_{k-1}^{(1)}, \ldots, X_{k-1}^{(N_p)} \overset{i.i.d.}{\sim} \pi_{k-1}^{N_p}$;  
        
		Sample $\widetilde{X}_k^{(v,1)}\hspace*{-0.4em}, \ldots, \widetilde{X}_k^{(v,N_p)} \overset{i.i.d.}{\sim}\hspace*{-0.1em} \text{P}^v_k(x_{k-1}^{(n)},\cdot),\ \forall v\in V$; 
		
		Compute $w_{k}^{(V_b,n)} = \frac{\prod_{v \in V_b} l^v\big( \tilde{x}^{(v,n)}_k,y^{(v)}_k \big)}{\sum_{n'=1}^{N_p} \prod_{v \in V_b} l^v \big( \tilde{x}^{(v,n')}_k,y^{(v)}_k \big)}$ for every $V_b \in \mathcal{K}$ and $n=1,\ldots,N_p$;

		Let $\pi_{k}^{N_p} = \textstyle\bigotimes_{V_b \in \mathcal{K}} \textstyle\sum_{n=1}^{N_p}w_k^{(V_b,n)}\delta_{\tilde{x}_k^{(V_b,n)}}$.
	}
\end{algorithm2e}

% ===============================================================================

%===============================================================================
%
%===============================================================================
\section{Case-study: The Oregonator system} \label{sec: Results}
%===============================================================================

In this section, we present the results obtained by the BPF (Section \ref{section: the block particle filter}) when used to reconstruct the concentrations from a benchmark reaction-diffusion system. Specifically, we focus on the task of estimating the state from a bistable, oscillatory Belousov-Zhabotinskii (BZ) reaction system \citep{Zhabotinsky1991}. 
The Oregonator \citep{Field1974} consists on the simplest realistic model of the BZ reaction dynamics, with network
\begin{align} \label{eq: Oregonator_CRN}%
    \begin{matrix}
        \begin{aligned}
            2\mathcal{S}_1                 &\,\xrightarrow{\kappa_1}\, \mathcal{S}_4 + \mathcal{S}_5\\
            \mathcal{S}_1 + \mathcal{S}_4  &\,\xrightarrow{\kappa_3}\, 2\mathcal{S}_1 + 2\mathcal{S}_2 \\
            \mathcal{S}_2 + \mathcal{S}_6  &\,\xrightarrow{\kappa_5}\, 0.5 \sigma \mathcal{S}_3,% ) <~ \sigma being positive is stated in the footnote. I Think it is better there than here
        \end{aligned}
        & &
        \begin{aligned}
            \mathcal{S}_1 + \mathcal{S}_3  &\,\xrightarrow{\kappa_2}\, 2\mathcal{S}_5 \\
            \mathcal{S}_3 + \mathcal{S}_4  &\,\xrightarrow{\kappa_4}\, \mathcal{S}_1 +  \mathcal{S}_5\\
            \phantom{\mathcal{S}_9 \, \xrightarrow{\kappa_6} \mathcal{S}_9} 
        \end{aligned}
    \end{matrix}
    % 
    % \nonumber\\
    % 
\end{align} 
with concentrations $z = \{ [\mathcal{S}_1],\ldots,[\mathcal{S}_6] \}$. Typically, species $(\mathcal{S}_4,\mathcal{S}_5,\mathcal{S}_6)$ are present in high densities and thus concentrations $z_{4,5,6}(s,t)$ are assumed constant on the timescale of a few oscillations. Moreover, when $z_3$ is slowly varying, the dynamics of the system are summarised by the evolution of $z_1$ and $z_2$ only, with $f_z$ represented in a scaled form\footnote[6]{From its mass-action dynamics, the Oregonator has dynamics $f_z(\cdot)$ in scaled form, defined by $[f_z]_1 = \epsilon^{-1} \big( z_1(1-z_1) - \frac{\sigma z_2(z_1-q)}{z_1+q} \big) + D_{z_1} \nabla^2 z_1$ and $[f_z]_2 = (z_1 - z_2) + D_{z_2} \nabla^2 z_2$
\citep{Keener1986}. 
The dimensionless constants are $(\epsilon, \sigma, q) = (0.08, 0.95, 0.0075)$. As for the diffusion coefficients, $(D_{z_1}, D_{z_2}) =  (5\times 10^{-4}, 5\times 10^{-6})$.
}.

We consider the system on a quasi-two-dimensional space $U$ with size $\overline{U} = 2$ units-of-space. We discretise $U$ with $\Delta u = 0.02$ units-of-size and consider the finite-difference approximation of the dynamics (Section \ref{sec: Reaction_Diffusion_FDM}) onto lattice $V$ with $\overline{V} = \overline{U}/\Delta u = 100$ grid-points in each dimension. The dynamics are then discretised in time with $\Delta t = 0.01$.
The process is assumed to be subjected to driving noise with coefficient
$
    g(X_{k-1}) =  \sigma_x \texttt{diag}(X_{k-1}) 
$ 
given $\sigma_x = 10^{-2}$.
The measurement process is defined by $H = [I_{\overline{V}^2}\otimes\Phi_{\mathcal{S}_1}\ I_{\overline{V}^2}\otimes\Phi_{\mathcal{S}_2}]$ with spectra $(\Phi_{\mathcal{S}_1},\Phi_{\mathcal{S}_2})$ collected at 10 equally-spaced wavelengths $\lambda \in [0,50)$ through response functions $\phi_{S_1}(\lambda) = \exp{\big({-}\frac{(\lambda - 10)^2}{30}\big)}$ and $\phi_{S_2}(\lambda) = \exp{\big({-}\frac{(\lambda - 40)^2}{30}\big)}$. The measurement noise is $e_k \sim \mathcal{N}(0,\sigma_y^2I_{N_y})$ with $\sigma_y^2 = 10^{-5}$. Finally, we assume initial condition $X_0 \sim \delta_{x_0}$, with $\delta$ denoting the Dirac delta distribution and $x_0$ is the non-trivial steady-state solution the Oregonator dynamics. 

For the estimation task, we design a BPF by partitioning the lattice $V$ into blocks with size $\overline{V}_b = 5$ (with a total of $(\overline{V}/\overline{V}_b)^2 = 400$ blocks). For each block, we consider an ensemble of $N_p = 128$ particles, each assumed to evolve according to the same dynamics as the simulation system. We analyse the filtering results when particles and importance weights are obtained using the conventional ($p(x_k|x_{k-1})$) or the optimal (Eq. \ref{eq: SIR-proposal}) importance density. 

\subsection{Estimation Results}

\begin{figure*} \centering
    \includegraphics[width=\linewidth]{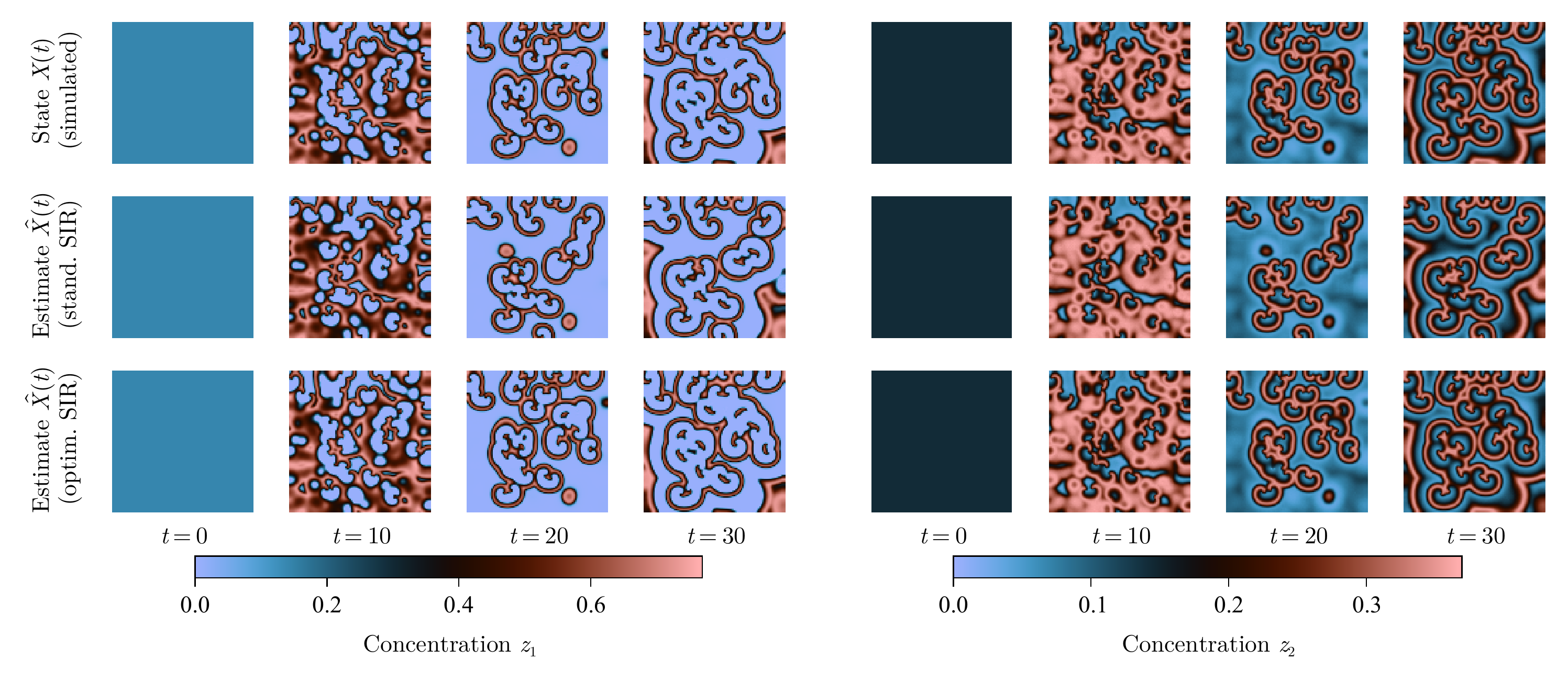}
    \vspace*{-2.75em}
    \caption{Snapshots of $[X_k]_1$  (4 first columns) and $[X_k]_2$ (4 last columns) for the Oregonator model at different times $t_k$. 
    First row shows a draw for $X_k$ used to generate observation $y_k$. 
    We estimated $\hat{X}_k$ by the BPF with a standard (second row) 
    and an optimal (third row) choice of importance sampling distribution with resampling steps (SIR).}
    \label{fig: estimation for the oregonator - standard SIR}
\end{figure*}

We apply the BPF with the described configuration to estimate concentrations $z = (z_1,z_2)$ from the Oregonator model using the generated spectral data. 
The results are shown in Fig.~\ref{fig: estimation for the oregonator - standard SIR} for different times $t_k$ through the mean of the approximated conditional distribution $\pi_k$, i.e. $\hat{X}_k^{(V_b)} := \sum_{n=1}^{N_p} w_k^{(V_b,n)}X_k^{(V_b, n)}$, and the actual state $X_k$.
From the evolution of the system, we observe the formation of oscillating ``spiral-patterns'' after a few iterations. The system is expected to switch back and forth between activator dominated states (species $\mathcal{S}_1$ is present in high-concentrations) and inhibitor dominated states (species $\mathcal{S}_2$ is present in high-concentrations). Initially, we observe concentrations to oscillate around the entire spatial domain. In this transient period ($t < 20$ units-of-time), the emerging patterns continuously change their shape and position. Afterwards, the patterns can be observed to emerge around fixed locations (as observed in $t \geq 20$ units-of-time). 
The results indicate that the BPF using the standard choice of importance distribution obtains poor estimates of the state: The particle approximations diverge during the transient region and ultimately form different patterns compared with those of the actual system. Conversely, the filter using the optimal importance distributions is able to accurately estimate the pattern-formation from the system. 

\begin{figure}[h!] \centering
    \includegraphics[width=\columnwidth]{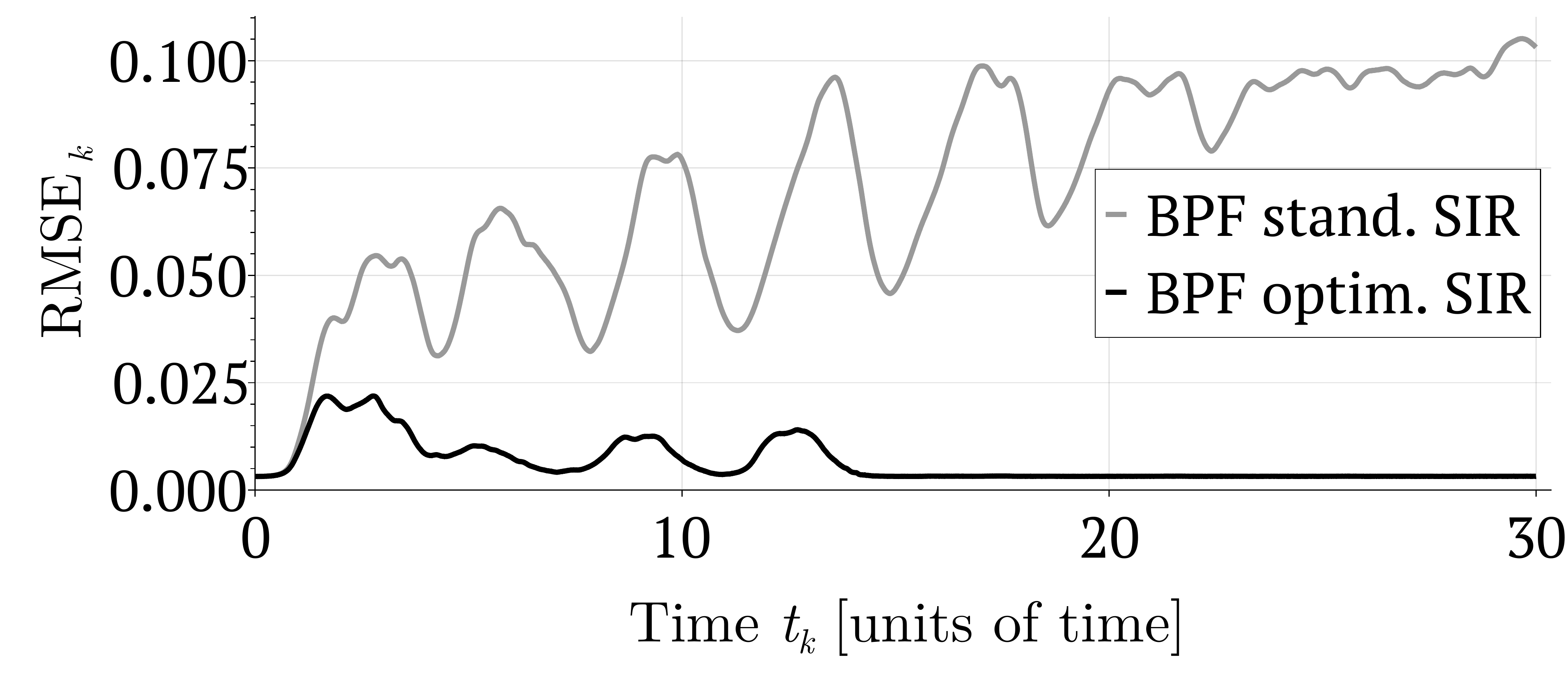}
    \vspace*{-2.1em}
    \caption{$\text{RMSE}_k$ between observation $y_k$ and prediction $\hat{y}_k$.}
    \label{fig: RMSE standard vs optimal SIR}
\end{figure}

In the following, we investigate the mismatch in Fig.~\ref{fig: estimation for the oregonator - standard SIR} by comparative metrics between both approximations. 
For each block $V_b \in \mathcal{K}$, we quantify the estimation accuracy in terms of the Root Mean Square Error (RMSE),
\begin{align}
    \text{RMSE}_k^{(V_b)} = \sqrt{\|\hat{y}_k^{(V_b)}{-}y_k^{(V_b)}\|^2_2},
\end{align}
with predicted outputs $\hat{y}_k^{(V_b)} := \sum_{n=1}^{N_p} w_k^{(V_b,n)}(Hx_k^{(V_b, n)})$, observation $y_k$, and $w_k^{(V_b, n)}$ and $x_k^{(V_b, n)}$ denoting, respectively, the normalised weight and the realisation of the $n$-th particle $X_k$ in a block $V_b$. Adding the accuracy obtained at each block, we arrive at the results shown in Fig.~\ref{fig: RMSE standard vs optimal SIR}. 

The results (Fig.~\ref{fig: RMSE standard vs optimal SIR}) show that the estimation accuracy oscillates due to the bistability of the system; crests relate to periods in which either the activator or inhibitor are dominated states, and troughs correspond to the transient periods. Moreover, this metric further highlights that accurate estimates are only obtained when the filter considers the optimal importance distribution. In this case, the estimation error increases slightly during $t \leq 15$, before decreasing to almost zero. Conversely, we observe the estimation accuracy to deteriorate quickly when the standard importance distribution is used.

\begin{figure}[h!] \centering
    \includegraphics[width=\columnwidth]{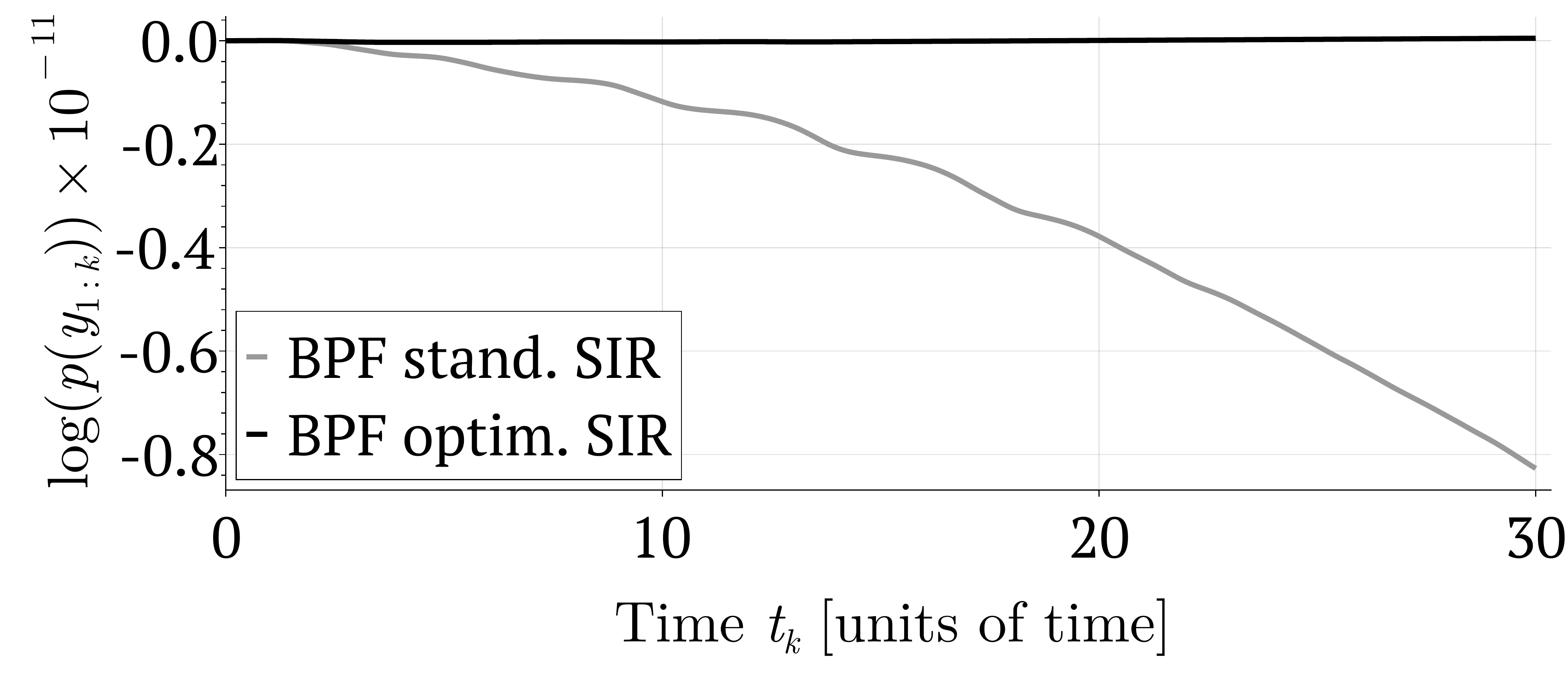}
    \vspace*{-2.1em}
    \caption{Estimated log-marginal likelihood $\texttt{log}(p(y_{1:k}))$. }
    \label{fig: p(y) standard vs optimal SIR}
\end{figure}

Furthermore, we compare the performance of both filters by computing the marginal likelihood function $p(y_{1:K})$ through Monte Carlo integration over the state space. Specifically, for each block $V_b\in\mathcal{K}$ we approximate the density
\begin{align}
p(y^{(V_b)}_{1:K}) \approx \prod_{k=1}^K \left(\cfrac{1}{N_p}\sum_{n=1}^{N_p}w^{(V_b,n)}_k\delta_{X^{(V_b,n)}_{k}}\right).
\end{align}
Avoiding numerical issues, we approximate the marginal likelihood in the logarithmic scale. Finally,  $\texttt{log}(p(y_{1:K}))$ over the full output-space is approximated by adding together the block-wise computed values, shown in Fig.~\ref{fig: p(y) standard vs optimal SIR}.

Here, the marginalised likelihood is the probability of the collection of observations $y_{1:K}$ given the model in Eq.~\eqref{eq: SS}, 
without assuming a particular realisation of the state. 
Under such definition, approximations of $p(y_{1:k})$ should have similar results regardless of which importance 
sampling distribution is adopted for the filtering, provided that enough particles are used to avoid degeneracy. However, the approximation obtained from the filter with the standard choice does not match that obtained from its optimal counterpart. 
The relatively small probability densities of the collection of observations in Fig.~\ref{fig: p(y) standard vs optimal SIR} follows as a direct consequence of path degeneracy.

%===============================================================================
\section{Concluding remarks} \label{sec: Concluding remarks}
%===============================================================================

In this work, we highlighted a sequential Monte Carlo approximation to the solution to the filtering equations~\eqref{eq: prediction} and \eqref{eq: correction}. Under the context of a reaction-diffusion system, we examined how the BPF  dealt with the high-dimensionality of the state under a standard and an optimal choice of importance sampling distribution.

BZ reaction systems are characterised by long-term unpredictability arising from an extreme sensitivity to initial conditions and process noise. Depending on these parameters, the system may or may not present oscillating patterns as in Fig.~\ref{fig: estimation for the oregonator - standard SIR}. Here, we presented the scenario of reconstructing the concentrations in the case of moderate values of the system noise variance. We would like to remark that the benefits from the optimal proposal will be negligible as this noise becomes small, as discussed in \cite{Snyder2015}. 

The BPF has some intrinsic limitations by design. For example, the bias introduced as a result of the blocking operator is not spatially homogeneous. In an extended version of this work, we intend to (i) assess any spatial inhomogeneity, (ii) investigate how some extensions to the BPF address the bias, and (iii) draw comparisons to other high-dimensional filters in a broader experimental study.  

%===============================================================================

\bibliography{IFAC_BPF4ReacDiff}
            
\end{document}